\documentclass[final,5p,times,twocolumn]{elsarticle}

\usepackage{amssymb}
\usepackage{amsmath}
\usepackage[english]{babel}
\usepackage{url}
\usepackage{xcolor}
\usepackage{paralist}

\newcommand{\contr}[1]{\mathcal{X}_{#1}}
\newcommand{\tra}[3]{\tilde{#1}_{#2:#3}}
\newcommand{\pare}{^\text{P}}

\newcommand{\chil}{^\text{C}}

\newcommand{\term}{^\text{f}}
\newcommand{\pmpc}{P-MPC}
\newcommand{\cmpc}{C-MPC}
\newtheorem{theorem}{Theorem}
\newtheorem{proof}{Proof}
\newtheorem{remark}{Remark}

\newcommand\norm[1]{\left\lVert#1\right\rVert}
\newcommand{\bx}{\boldsymbol{x}}
\newcommand{\bu}{\boldsymbol{u}}
\newcommand{\np}{N}
\newcommand{\Pnp}{N^{\text{P}}}
\newcommand{\Cnp}{N^{\text{C}}}
\newcommand{\wP}{\boldsymbol{w}^{\text{P}}}
\newcommand{\uP}{\boldsymbol{v}^{\text{P}}}
\newcommand{\xP}{\boldsymbol{z}^{\text{P}}}
\newcommand{\termPar}{\mathcal{Z}^{\text{f,P}}}
\newcommand{\eP}{\boldsymbol{e}^{\text{P}}}

\newcommand{\EP}{\mathcal{E}^{\text{P}}}
\newcommand{\EC}{\mathcal{E}^{\text{C}}}
\newcommand{\uC}{\boldsymbol{v}^{\text{C}}}
\newcommand{\xC}{\boldsymbol{z}^{\text{C}}}
\newcommand{\Kf}{K^{\text{f,P}}}
\newcommand{\Vc}{\mathcal{V}^{\text{C}}}
\newcommand{\KP}{K^{\text{P}}}
\newcommand{\Ts}{T^{\text{s}}}

\allowdisplaybreaks
\journal{Systems \& Control Letters}

\begin{document}

\begin{frontmatter}



\title{Recursive Feasibility without Terminal Constraints via Parent-Child MPC Architecture}


 \author{Filip Surma\corref{cor1}\fnref{delftlabel}}
 \author{Anahita Jamshidnejad\fnref{delftlabel}}
\cortext[cor1]{Corresponding author.\\ E-mail addresses: f.surma@tudelft.nl (F. Surma), a.jamshidnejad@tudelft.nl (A. Jamshidnejad).}


\affiliation[delftlabel]{organization={Delft University of Technology},
            addressline={Kluyverweg 1}, 
            city={Delft},
            postcode={2629 HS}, 
            state={South Holland},
            country={The Netherlands}}

\begin{abstract}
This paper proposes a novel hierarchical model predictive control (MPC) framework, 
called the Parent-Child MPC architecture, to steer nonlinear systems 
under uncertainty towards a target set, balancing computational complexity and guaranteeing recursive 
feasibility and stability without relying on conservative terminal constraints in online decision-making. 
By coupling a small-horizon Child MPC layer with one or more large-horizon Parent MPC layers, 
the architecture ensures recursive feasibility and stability through adjustable stage-wise constraints 
derived from tube-based control. 
As is demonstrated in our case studies, compared to traditional MPC methods, 
the proposed Parent-Child MPC architecture enhances performance and computational efficiency, reduces conservativeness, and 
enables scalable planning for certain nonlinear systems.
\end{abstract}


\end{frontmatter}



\section{Introduction}
\label{sec:intro}

Model predictive control (MPC) employs a model of the system to predict its future 
behavior over a finite time horizon and to accordingly optimize the control input trajectory. 
A key strength of MPC lies in its ability to explicitly incorporate constraints during the optimization procedure \cite{MPCbook}. 

As infinite constrained optimization is usually computationally intractable, 
MPC  operates over a finite prediction horizon. 
Due to this, ensuring stability for MPC is generally challenging, and 
is commonly achieved via, e.g., enforcing passivity conditions  \cite{passive},   
using the augmented stage cost function corresponding to a single time step as a control Lyapunov function \cite{stateCostModifided},  
or treating the optimal cost of the MPC optimization problem as a control Lyapunov function and  
demonstrating its decline \cite{MPCbook}.%

These methods for MPC often involve supplementary, terminal constraints. 
Applied exclusively to the final predicted state, terminal constraints guide the system towards a terminal state region 
$\mathcal{X}^{\text{f}} \subseteq \mathcal{X}$ near its destination (with $\mathcal{X}$ 
the admissible state set). 
The $\np$-step controllable set $\contr{\np}$ set is a set including all states  
that are steerable by an admissible control input sequence of length $\np$ 
---  the number of time steps across the prediction horizon --- into 
 set $\mathcal{X}^{\text{f}}$ \cite{Invariant}. 
For larger prediction horizons, this controllable set is usually expanded. In other words, 
some states, although safely steerable into $\mathcal{X}^{\text{f}}$ with large prediction horizons,  
lie outside the reach of short-horizon MPC.%

Large-horizon MPC is typically computationally prohibitive \cite{extendedHorizon}.  
It is thus common to simplify (e.g., linearize) the prediction model \cite{quadrotor,mobileRobot}. 
This, however, potentially leads to constraint violations, due to model inaccuracies.
Alternatively, as explained in \cite{sampling}, the sampling time is manipulated  
to enlarge the controllable state set without increasing the prediction horizon.%

Hybrid approaches can integrate MPC with other controllers. 
An example involves dividing the prediction horizon into a control horizon governed by MPC  
and a segment where the last MPC input is simply repeated \cite{overviewMPC}. 
This approach usually enhances computational efficiency by reducing the number of decision variables. 
Alternatively, after executing the MPC trajectory 
for given time steps, a terminal control law (e.g., a linear quadratic regulator  
for linear systems) is activated \cite{MPCbook,mobileRobot,DualTerminal}. 
This method is based on invariance of the terminal set.%

These strategies all rely on reaching a given terminal set via MPC, without 
inherently enlarging the $\np$-step controllable state set. Hierarchical MPC offers promising alternatives: 
In \cite{Hier,Hier2}, MPC performs high-level planning, generating optimal reference trajectories.  
A lower-level controller executes the plan. 
A recent bi-level MPC architecture in \cite{Anahita}   
synergizes two MPC formulations: ``Rough long-vision'' MPC with a coarse prediction model and large sampling time, 
and ``detailed short-vision'' MPC for enhancing and executing the plan across a smaller prediction horizon with smaller sampling times. 
Long-term MPC provides a trajectory that guides the constraints of short-term MPC. 
This architecture has been extended to long-term exploration planning for multi-robot systems in \cite{FLMPC}, 
with highly promising results, but no formal stability guarantees.%

Inspired by this, we address the challenge of reaching a designated terminal set through online 
decision-making by proposing a multi-level architecture composed of interconnected MPC systems, 
leveraging robust Tube-based MPC (TMPC). 
TMPC ensures constraint satisfaction under bounded uncertainties  
without significantly increasing online computational complexity \cite{MPCbook}. 
It uses a fixed-shape tube (e.g., a polyhedron or ellipsoid
\cite{setComputation}) around the nominal trajectory to guarantee safety \cite{LTMPC}. 
The nominal component of TMPC governs the center of the tube, while an   
ancillary controller maintains the actual state within the tube. 
The fixed tube structure, however, may lead to conservative behavior and reduced performance. 
For instance, applied to a wheeled E-puck robot \cite{epuck}, 
TMPC caused a substantial unnecessary reduction in the robot velocity, up to $30\%$ \cite{SDDTMPC}.  
This performance degradation stems from the dual role of the ancillary control law, 
as it should both ensure stability and compensate for disturbances and model discrepancies. 
These dual objectives restrict the flexibility of nominal MPC in selecting optimal control inputs. 
Our key contributions paper are: 
\begin{compactitem}
    \item   
We introduce a   novel hierarchical control framework, the  
Parent-Child MPC architecture, which effectively replaces the nominal component in 
TMPC for certain nonlinear systems. It guarantees convergence to  
a designated terminal set without imposing a terminal constraint, 
thereby avoiding conservativeness and risk of infeasibility typically associated with such constraints. 
\item  
A small-horizon Child MPC (\cmpc) controller directly steers the system, emphasizing local optimality. 
One or more large-horizon Parent MPC (\pmpc) layers address stability and long-term feasibility. 
Strategic interactions of these layers allow to replace terminal constraints with 
adjustable stage-wise state constraints, including robust positive invariant tubes 
derived by Parent layers. 
This significantly enlarges the 
feasible exploration space of \cmpc, leading to improved performance. 
\item  
The cross section of these tubes is composed of a fixed subset $\EP$ of the state space. We prove that $\EP$ 
is also a robust positive invariant set for the C-MPC nominal trajectory. 
Once this nominal trajectory enters $\EP$, the goal of Parent-Child MPC architecture is met, 
allowing a secondary low-complexity controller (e.g., small-horizon conventional TMPC) 
to steer the actual state towards a neighborhood of the origin.
\end{compactitem}
Similarly to \cite{Anahita}, the trajectories of \pmpc\ layers --- based on linear TMPC \cite{LTMPC} ---  
are used to formulate constraints for \cmpc\  
that ensure long-term feasibility and stability. 
\cmpc, similarly to an ancillary controller, directly steers the system, 
which is initialized far from its desired state, with no reference trajectory. 
Unlike the approach in \cite{NTMPC}, where MPC serves as an ancillary controller 
within a nonlinear TMPC framework, \cmpc\  
does not track the nominal trajectories of \pmpc.  
Instead, it shares the same cost function as \pmpc, while 
employing a more accurate prediction model and 
finer sampling times.%

Section~\ref{sec:fromulation} formulates the problem. Section~\ref{sec:Method} details the 
Parent-Child MPC architecture. 
Section~\ref{sec:implemenation} presents and discusses the results of a case study, 
and Section~\ref{sec:conclusions} concludes the paper and proposes directions for 
future research.%

\section{Problem formulation}
\label{sec:fromulation}

This section describes the problem of controlling a nonlinear system subject to additive disturbances, 
with the goal of steering the system from an arbitrary initial state to a desired target state. 
The system model follows the formulation in \cite{LipsitzImplmented,Lipschitz}.%

\subsection{System modeling}
\label{sec:sys_model}

The system dynamics is described in discrete time by the following nonlinear model subject to additive disturbances:
\begin{equation}\label{eq:system}
    \boldsymbol{x}_{k+1}=A\boldsymbol{x}_k +g(\boldsymbol{x}_k)+B\boldsymbol{u}_k+\boldsymbol{w}_k
\end{equation}
where $\boldsymbol{x}_k \in \mathbb{R}^n \subseteq \mathcal{X}$ is the state vector, 
$\boldsymbol{u}_k \in \mathbb{R}^m \subseteq \mathcal{U}$ is the control input vector, 
$\boldsymbol{w}_k \in \mathbb{R}^n \subseteq \mathcal{W}$ is an unknown bounded disturbance vector. 
The states, control inputs, and disturbances are constrained to convex, compact sets $\mathcal{X}$, $\mathcal{U}$, and $\mathcal{W}$, respectively.  
Function $g(\cdot)$ is nonlinear Lipschitz with Lipschitz constant $L>0$, satisfying $g(\boldsymbol{0})=\boldsymbol{0}$.  
The matrix pair $(A,B)$ is stabilizable. Lipschitz continuity of $g(\cdot)$ implies that:
\begin{equation}\label{eq:Lipschitz}
    \norm{g(\boldsymbol{x}_1)-g(\boldsymbol{x}_2)}\leq L \norm{\boldsymbol{x}_1-\boldsymbol{x}_2}
    \qquad  \forall \boldsymbol{x}_1, \boldsymbol{x}_2\in\mathcal{X}    
\end{equation}

\subsection{Control problem}
 
The control objective is to regulate the state towards 
a suitable robust positive invariant set including 
the origin, as defined in \cite{LTMPC}. 
Until this goal is reached, a stage cost $l(\cdot,\cdot)$, given below, is accumulated, 
with $Q$ and $R$ positive definite matrices:
\begin{equation}\label{eq:cost}
    l(\boldsymbol{x}_k,\boldsymbol{u}_k)=\boldsymbol{x}_k^\top Q\boldsymbol{x}_k+\boldsymbol{u}_k^\top R\boldsymbol{u}_k
\end{equation}
As shown in \cite{Lipschitz}, systems given by \eqref{eq:system} under the stated conditions, 
can robustly be controlled using TMPC. 
The nominal optimization problem of nonlinear TMPC at time step $k$, with prediction horizon $\np$, for system \eqref{eq:system} 
is given by:
\begin{subequations}
    \label{eq:TMPC}
    \begin{align}  
    \begin{split}\label{cost_fun_TMPC}
        V(\boldsymbol{x}_k)=\min_{\tra{\boldsymbol{z}}{k}{k+\np},\tra{\boldsymbol{v}}{k}{k+\np-1}}\sum^{k+\np - 1}_{i=k}
        \left(\boldsymbol{z}^\top_{i|k}Q\, \boldsymbol{z}_{i|k}+\boldsymbol{v}^\top_{i|k}R\, \boldsymbol{v}_{i|k}\right)+
        \boldsymbol{z}^\top_{k+\np|k}P\boldsymbol{z}_{k+\np|k} 
    \end{split}
    \\
    &\text{such that for $i=k,...,k + \np-1$:} 
\nonumber
\\
    \begin{split}\label{initi_TMPC}
        \boldsymbol{{x}}_k\in \{\boldsymbol{z}_k\} \oplus \mathcal{E}
    \end{split}
    \\
    \begin{split}\label{model_TMPC}
        \boldsymbol{z}_{i+1 | k}=A\boldsymbol{z}_{i|k} +g\left(\boldsymbol{z}_{i|k}\right)+B\boldsymbol{v}_{i|k} 
    \end{split}
    \\
    \begin{split}\label{state_TMPC}
        {\boldsymbol{z}}_{i+1|k}\in \mathcal{X}\ominus \mathcal{E}
    \end{split}
    \\
    \begin{split}\label{input_TMPC}
        {\boldsymbol{v}}_{i|k}\in \mathcal{V}
    \end{split}
    \\
    \begin{split}\label{ter_TMPC}
        {\boldsymbol{z}}_{ k + \np | k}\in \mathcal{Z}\term
    \end{split}
\end{align}
\end{subequations}
which, rather than directly steering the system state $\boldsymbol{x}_k$, 
determines nominal control input $\boldsymbol{v}_k$ 
to steer the nominal state $\boldsymbol{z}_k$ . 
The nominal prediction model is given via \eqref{model_TMPC},  
which eliminates the disturbance term $\boldsymbol{w}_k$ from \eqref{eq:system}. 
The nominal state trajectory represents the center of a rigid tube with cross section    
$\mathcal{E}$, i.e., the set of possible errors between nominal and actual states. 
The state admissibility constraint \eqref{state_TMPC} has been tightened 
to ensure the actual state, impacted by uncertainties, remains within the admissible state set $\mathcal{X}$. 
The notation $a_{k_2|k_1}$, with $k_2 \geq k_1$, shows the value of variable $a$ at time step $k_2$  
predicted at time step $k_1$, with $a_{k_1|k_1} : = a_{k_1}$.  
Tilde for decision variables  
shows a sequence, where $\tra{a}{k_3}{k_4}$ with $k_4 \geq k_3$ represents $[a_{k_3},\dots,a_{k_4}]^\top$.  
Constraint \eqref{initi_TMPC} adjusts the nominal state $\boldsymbol{z}_{k}$ at time step $k$  
to locate the cross section $\{\boldsymbol{z}_k\} \oplus \mathcal{E}$ of the tube 
within set $\mathcal{X}$, so that the measured state $\bx_k$ lies within the tube. 
For all time steps $i=k,...,k + \np-1$, the actual control input $\boldsymbol{u}_{i|k}$ should satisfy $\boldsymbol{u}_{i|k} \in \mathcal{U}$. 
This requires defining tightened constraints on nominal control inputs, as is done in \eqref{input_TMPC}. 
The terminal constraint \eqref{ter_TMPC}, where $\mathcal{Z}\term$ is a control invariant set, designed as 
a tightened version of the actual terminal state set $\mathcal{X}^{\text{f}}$,  
together with the terminal cost $\boldsymbol{z}^\top_{k + \np | k}P\; 
\boldsymbol{z}_{k + \np | k}$ in \eqref{cost_fun_TMPC}, 
ensures recursive feasibility --- i.e., if a solution initially exists, a solution exists for all future time steps.  
However, as explained earlier, imposing such a terminal constraint may restrict the set of initial states 
from which the system is steerable to the origin. 
Note that $\oplus$ and $\ominus$ show the Minkowski set summation and subtraction, respectively.%

The actual system, which is prone to disturbances and model discrepancies, is controlled via 
an ancillary control law:
\begin{align}
\label{eq:actual_control_input}
\boldsymbol{u}_k=\pi(\boldsymbol{x}_k,\boldsymbol{z}_k,\boldsymbol{v}_k),
\end{align}
designed to ensure, through the tightened constraint \eqref{state_TMPC}, that the evolved actual state $\boldsymbol{x}_{k+1}$ 
remains within the tube, provided that the current state $\boldsymbol{{x}}_k$ satisfies \eqref{initi_TMPC}. 
The tube cross section $\left\{\boldsymbol{z}_{i|k} \right\} \oplus \mathcal{E}$, for all time steps $i=k,...,k + \np-1$, 
is a robust positive invariant set when the system follows control input $\pi\left(\boldsymbol{x}_k,\boldsymbol{z}_k, \boldsymbol{v}_k \right)$.
 
Our goal is to design an MPC-based control architecture that can steer the system to the origin 
from any admissible initial state --- assuming a feasible control problem --- while guaranteeing stability and recursive feasibility. 
Although solving \eqref{eq:TMPC} with a sufficiently large prediction horizon theoretically achieves this goal,  
the computational burden may be prohibitive \cite{MPCbook,SDDTMPC,FLMPC}. The control architecture proposed in this paper 
addresses this issue by alternately solving an additional quadratic program alongside the nonlinear program, 
reducing the computational complexity while maintaining the desired properties.%

\section{Proposed architecture: Parent-Child MPC}\label{sec:Method}

This section discusses the proposed Parent-Child MPC architecture, 
which replaces the nominal problem in conventional TMPC given by \eqref{eq:TMPC}  
to offer a computationally efficient alternative with stability guarantees. 
The Parent MPC (\pmpc) extends standard linear TMPC, incorporating both a nominal 
MPC problem and a feedback ancillary control law. 
In contrast, the Child MPC (\cmpc) leverages the nominal problem from the nonlinear 
TMPC formulation in \eqref{eq:TMPC} to more flexibly steer the system. 
Unlike traditional TMPC, \cmpc\ omits terminal constraints, 
thereby eliminating the restrictive influence of these sets, as discussed in Section~\ref{sec:intro}. 
Instead, \cmpc\ operates under additional stage-wise constraints derived by \pmpc, 
in the form of translated tubes -- robust positive invariant sets -- 
composed of a fixed subset $\EP$ of the state space centered along the Parent nominal state trajectory.%

Set $\EP$ is designed to lie within the controllable state set, from which 
the nominal Child trajectory is steerable into a desired terminal set.  
We prove that, under the Parent-Child MPC strategy, the untranslated set $\EP$ itself is 
a robust positive invariant set for the nominal \cmpc\ trajectory. 
Therefore, once the nominal \cmpc\ trajectory enters this set, 
the primary objective of the Parent-Child MPC architecture is fulfilled. 
From this point onward, a secondary control phase -- such as a conventional TMPC scheme 
with a small horizon (and thus affordable online computations), 
or another control-invariant strategy -- may be employed to steer the actual state trajectory 
towards a desired neighborhood of the origin.%

\pmpc\ employs a simplified nominal version of the system model 
in \eqref{eq:system} that eliminates both the Lipschitz nonlinearity and the additive disturbances. 
The underlying linear TMPC formulation of \pmpc\ is designed 
to mitigate the uncertainties introduced by the model discrepancy 
that stems from omitting the Lipschitz nonlinearities of the original system. 
By adopting a coarser sampling time for prediction and optimization,  
\pmpc\ reduces computational complexity and supports larger prediction horizons. 
While this long-term planning facilitates stability, directly following the trajectories 
generated by \pmpc\ may lead to suboptimal performance due to model inaccuracies.%
  
\cmpc\ operates at a finer resolution, using a nominal, yet nonlinear, version of 
the model in \eqref{eq:system} that solely eliminates the additive disturbances. 
Exploiting its more accurate model, \cmpc\ focuses on short-term optimality 
within a disturbance-driven tube, 
while inheriting stability and feasibility guarantees from \pmpc. 
An ancillary control strategy (such as the one in \cite{NTMPC}) is ultimately employed to reject these disturbances.%

The Parent-Child MPC architecture integrates long-term planning capabilities of \pmpc\ 
with responsiveness and short-term optimality of \cmpc, 
while preserving stability and recursive feasibility. 
The primary distinction between conventional nonlinear TMPC 
and the Parent-Child MPC architecture lies in the role of the nominal components: 
In standard TMPC an ancillary controller ensures stability 
by keeping the state within a tube located in the state space 
through the nominal controller. 
Thus, reducing the tube size is essential to enhancing flexibility in the nominal trajectory, 
thereby improving the performance 
(see tightened state admissibility constraint \eqref{state_TMPC}). 
In contrast, within the Parent-Child MPC architecture, 
\cmpc\ enhances the optimality, 
whereas \pmpc\ ensures stability and recursive feasibility by generating a 
suitable modeling error tube. 
Consequently, enhancing the flexibility of \cmpc\ and, thus, the system performance, 
requires an enlarged modeling error tube.%

Next, we present the formulation of \pmpc\ and \cmpc. 
The actual state at time step $k$ is denoted by $\bx_k$, 
which is guided by the actual control input $\bu_k$. 
The nominal state and control input of \pmpc\ and \cmpc\ 
are represented by $\xP_k,\uP_k$ and $\xC_k,\uC_k$, respectively. 
The prediction horizon and cross section of the tube --- 
a robust positive invariant set for rejecting uncertainties --- for \pmpc\ and 
\cmpc\ are $\Pnp$, $\EP$ and $\Cnp$, $\EC$, respectively.  
The actual state and input of \pmpc\ are given by $\boldsymbol{x}_k\pare$ and $\boldsymbol{u}_k\pare$.%

\subsection{Parent MPC}
\label{sec:pmpc}

The formulation of \pmpc\ is based on standard linear TMPC \cite{LTMPC}. 
Rejecting the additive disturbances 
$\boldsymbol{w}_k$ affecting the original system modeled by \eqref{eq:system} is delegated to \cmpc. 
Consequently, external disturbances do not appear in the formulation of \pmpc, 
which is designed for the following dynamic model: 
\begin{align}
\label{eq:parent_model}
    \boldsymbol{x}\pare_{k+1} = A\boldsymbol{x}\pare_k + B \boldsymbol{u}\pare_k + \wP_k   
\end{align}
with $\wP_k = g\left(\boldsymbol{x}\pare_k\right) $, i.e., both the additive disturbance $\boldsymbol{w}_k$  
and the nonlinear term $g\left(\boldsymbol{x}\pare_k\right)$ as a controllable part of the dynamics are eliminated. 
Instead, the nonlinear term is treated as a bounded additive disturbance $\wP_k$. 
Since $g(\cdot)$ is Lipschitz continuous and $g(\boldsymbol{0}) = \boldsymbol{0}$, 
we have $\norm{ g\left(\boldsymbol{x}\pare_k\right)} \leq L\boldsymbol{x}\pare_k$. 
Given the boundedness of the state space $\mathcal{X}$, this ensures that $g\left(\boldsymbol{x}\pare_k\right)$ 
is also bounded.

The model in \eqref{eq:parent_model}, therefore, structurally 
resembles a linear system subject to bounded additive disturbances,  
consistent with the framework in \cite{LTMPC}. 
The nominal component of \pmpc\ then follows a disturbance-free version of \eqref{eq:parent_model}, given by:
\begin{equation}
\label{eq:simplfied}
    \xP_{k+1}=A \xP_k+B \uP_k
\end{equation}
In this context, the error $\eP_k$ between the actual model \eqref{eq:parent_model} 
and the nominal model \eqref{eq:simplfied} of \pmpc\ is defined as: 
\begin{align}
\label{eq:modeling_error}
\eP_k:=\boldsymbol{x}\pare_k - \xP_k
\end{align}

The simplified deign of \pmpc\ enables its nominal optimization problem to be cast as 
a convex or quadratic program, significantly reducing computational complexity.%

In linear TMPC,  the nominal control input is augmented 
with an ancillary control law in the form of linear state error feedback  
for regulating the modeling error $\eP_k$. Accordingly, the control input in \pmpc\ is defined as:
\begin{equation}
\label{eq:ancillary}
    \boldsymbol{u}_k\pare = \uP_k + \KP \eP_k
\end{equation}
where $\KP$ is a feedback gain matrix selected such that $A+B \KP$ is stable 
(i.e., all eigenvalues lie strictly within the unit circle). In general, matrix $\KP$ is not unique.%

The control law in \eqref{eq:ancillary} is not directly 
applied to the actual system. Instead, as is detailed in Section~\ref{sec:C-MPC}, it  
 serves as a baseline input around which the nominal \cmpc\ policy is constructed, 
 and warm-starts the \cmpc\ nominal optimization. 
Moreover, at the \pmpc\ level, the choice of the feedback gain matrix $\KP$ 
is primarily important for the construction of the    
modeling error tube with cross section $\EP$. This tube is subsequently injected into \cmpc\ 
as a sequence of stage-wise constraints on its nominal state, omitting the need for explicit terminal sets. 
This tube additionally plays a central role in establishing the recursive 
feasibility for the Parent-Child MPC architecture.%



\subsubsection{Nominal problem of \pmpc}
\label{sec:nominal_PMPC}

The nominal optimization problem of \pmpc\ spans the Parent prediction horizon $\Pnp$, 
and incorporates both the Parent tube with cross section $\EP$ and the Child tube 
with cross section $\EC$. 
This optimization problem at time step $k$ is given by:
\begin{subequations}\label{eq:Parent}
\begin{align}
    \begin{split}\label{cost_fun}
       & V\pare\left({\boldsymbol{z}}\chil_{k|k-1}\right)\\
        &=\min_{\tra{\boldsymbol{z}}{k}{k+\Pnp}\pare,\; \tra{\boldsymbol{v}}{k}{k+\Pnp-1}\pare} \Bigg\{\sum^{ k + \Pnp -1}_{i=k}
        \left({\xP_{i|k}} ^\top Q\, {\xP_{i|k}}+{\uP_{i|k}}^\top R\, {\uP_{i|k}}\right)\\
        & \hspace{30ex} + {\xP_{k + \Pnp|k}}^{\scriptstyle{\top}} P\, {\xP_{k + \Pnp|k} }\Bigg\}
    \end{split}
    \\
        &\text{such that for $i=k,...,k + \Pnp-1$:}
        \nonumber
        \\
    \begin{split}\label{initi}
      {\boldsymbol{z}}\chil_{k|k-1} \in \left\{\xP_{k} \right\}  \oplus \EP
    \end{split}
    \\
    \begin{split}\label{model}
        \xP_{i+1 | k} = A \xP_{i | k}+B \uP_{i | k} 
    \end{split}
    \\
    \begin{split}\label{state}
        \xP_{i+1|k}\in \left(\mathcal{X}\ominus\EC\right) \ominus \EP 
    \end{split}
    \\
    \begin{split}\label{input}
        \uP_{i|k}\in \Vc \ominus \KP \EP
    \end{split}
    \\
    \begin{split}\label{ter}
        \xP_{ k + \Pnp|k}\in \termPar
    \end{split}
    \end{align}
\end{subequations}
This problem formulation enhances flexibility by encoding the initial 
Parent nominal state $\xP_{k}$ as a decision variable. As enforced by constraint \eqref{initi},   
this state should be determined such that the most recent estimate  ${\boldsymbol{z}}^C_{k|k-1}$ 
of the Child nominal state for time step $k$ lies within the current Parent tube. 
Since the Parent tube is by design a robust positive invariant set, all future Parent nominal states 
are then guaranteed to remain within this tube. 
Note that the nominal Child state $\xC_{k}$ is a decision variable in \cmpc\ (see Section~\ref{sec:C-MPC}), 
which is solved only after the \pmpc\ optimization problem. 
This explains deploying the most recent available estimation $\xC_{k|k - 1}$  
as a proxy for the Child nominal state at time step $k$ in constraint \eqref{eq:Parent}.%

The evolution of the nominal state in \pmpc\ is governed by the simplified dynamics \eqref{eq:simplfied}, restated in constraint \eqref{model}.%

To ensure constraint satisfaction despite model simplifications, constraints    
\eqref{state} and \eqref{input} impose tightened bounds on state and control input trajectories, respectively. 
While the trajectories generated by \pmpc\ do not directly steer the actual system, 
they indirectly influence the actual trajectories by steering 
the nominal state trajectory of \cmpc\ (i.e., the center of the Child tube) 
to ensure it remains stage-wise within the Parent tube.   
This condition, as detailed in Section~\ref{sec:C-MPC}, is essential to eliminating 
terminal constraints from the \cmpc\ formulation.  
Moreover, \cmpc\ relies on a simplified model of the system that 
excludes the additive disturbances. 
This implies that constraint tightening is necessary to always guarantee state admissibility.  
Accordingly, as encoded in constraint \eqref{state}, the admissible state set $\mathcal{X}$ 
is first shrunk by the shape of set $\EC$ --- cross section of the disturbance-driven 
tube of \cmpc --- and is further eroded by set $\EP$ --- cross section of the modeling error tube. 
This ensures that the set $\left\{ \xP_{i|k} \right\} \oplus \EP \oplus \EC$, 
which contains the actual state trajectory for all $i = k,...,k + \Pnp-1$, remains a 
subset of admissible set $\mathcal{X}$. 
A similar tightening is applied to the nominal control inputs of \pmpc. 
In constraint \eqref{input}, the admissible control input set 
$\mathcal{U}$ is replaced by a tightened set $\Vc$, designed based on the actual 
control input policy of \cmpc\ for admissibility of its nominal control inputs 
(see \cite{NTMPC} for details), and is further eroded by the set $\KP \EP$, 
thereby reserving room for augmenting the Parent ancillary control law later.%

To ensure stability, the terminal state $\xP_{k+N^p|k}$ in \pmpc\ must lie within 
a terminal set $\termPar$, enforced via constraint \eqref{ter}. 
The associated terminal cost ${\xP_{k + \Pnp|k}}^\top P\, \xP_{k + \Pnp|k}$ is defined in quadratic form, 
with matrix $P$ a constant, positive definite matrix obtained by solving the 
following discrete-time Lyapunov equation: 
\begin{equation}\label{eq:Lapunov}
    \left(A+B \Kf\right)^\top P \left(A+B \Kf \right)+Q+{\Kf}^{\scriptstyle{\top}} R \Kf=P
\end{equation}
Matrices $A, B, \Kf$ define the closed-loop dynamics under a terminal control law 
that follows a state feedback policy with gain $\Kf$, while 
$Q$ and $R$ are the stage cost matrices. 
Note that the feedback gains $\KP$, used in \eqref{eq:ancillary} for the ancillary control law of \pmpc,  
and $\Kf$, used to determine the terminal cost, both need to stabilize 
the simplified nominal model \eqref{eq:simplfied}. 
However, they do not necessarily need to be identical. Their independent design 
can introduce flexibility, as gain $\KP$ influences the shape and size of $\EP$, i.e.,  
the modeling error set, which constrains the 
nominal trajectory of \cmpc, while gain $\Kf$ affects the size of 
the terminal set $\termPar$ for \pmpc, which determines the admissible endpoint of the 
nominal trajectory of \pmpc. 
This terminal set may be defined as a level set of the terminal cost, 
often taking the form of a polyhedron or an ellipsoid, 
and obtained by imposing ${\xP_{k + \Pnp|k}}^\top P\, \xP_{k + \Pnp|k} <c$, where $c > 0$ 
is a scalar that determines the size of the terminal set.%


\subsubsection{Designing the modeling error tube cross section $\EP$} 
\label{sec:PMPC_tube}

The modeling error tube in \pmpc\ should be  
a robust positive invariant set with cross section $\EP$, such that 
$\eP_k\in\EP$ implies $\eP_i\in\EP$, 
for $i > k$. 
Various state-of-the-art algorithms, e.g., in \cite{setComputation,tube}, 
are available for designing cross section of such sets, taking into account that this set 
is not necessarily unique for a given system subject to given disturbance bounds.%

In conventional TMPC, determining optimal shape and size for the tube entails a fundamental trade-off. 
Smaller tubes, on the one hand, increase flexibility in designing the nominal trajectory 
by enlarging the tightened state admissible set (see  
constraints \eqref{state_TMPC} and \eqref{state}). 
On the other hand, a smaller tube often necessitates a larger feedback gain matrix $\KP$ to ensure robust stability. 
This, in turn, tightens the admissible control input space (see \eqref{input}). 
Thus, either design choice may lead to an overall degradation of closed-loop performance.%

Since the tube of \pmpc\ encapsulates the nominal state trajectory of \cmpc, 
shrinking this tube --- while expanding the admissible state space for the nominal 
state trajectory of \pmpc\ --- has the opposite effect on the nominal state trajectory of \cmpc. 
The trajectory of \cmpc, however, should enhance optimality. 
Hence, increased flexibility in state space exploration is highly desirable. 
Accordingly, the Parent-Child MPC architecture pursues a fundamentally different design objective 
than that of conventional TMPC:  
It enlarges the cross section $\EP$ of the modeling error tube, 
provided that the computational demand remains tractable for solving the online optimization of \cmpc. 
Since, by the end of the primary control phase, the nominal \cmpc\ trajectory lies within set $\EP$, 
this set becomes the feasible exploration region for the conventional TMPC in the secondary phase of control. 
Therefore, the design trade-offs for set $\EP$ should also account for the computational demands of the conventional TMPC. 
In practice, for many systems described by \eqref{eq:system}, with reasonably simple terminal sets, 
if \cmpc\ (which operates under stricter constraints) is computationally affordable over horizon $\Cnp$, then so is the conventional TMPC.

The enlarged tube of \pmpc\ serves as a stage-wise constraint for \cmpc, 
thereby expanding its feasible region and enabling it to pursue performance-enhancing control strategies. 
To preserve the stability and recursive feasibility guarantees of the overall framework, the set $\EP$ should satisfy:
\begin{align}
    \label{eq:EP_stability}
    \EP \subseteq\mathcal{Z}_{\Cnp}^{\text{C}}
\end{align}
where $\mathcal{Z}_{\Cnp}^{\text{C}}$ is the $\Cnp$-step controllable set for the nominal component of \cmpc\  
that steers its nominal state -- within $\Cnp$ time steps -- into a terminal region $\mathcal{Z}^{\text{f,C}}$, which is designed offline by tightening the actual 
terminal state set $\mathcal{X}^{\text{f}}$ to account for the impact of the additive 
disturbances on the actual system.%

This approach draws inspiration from applications of TMPC in passenger car systems \cite{LipsitzImplmented} and 
planar quadrotor drones \cite{LipsitzDownsied}.%

\begin{remark}
When computing the maximal set $\mathcal{Z}_{\Cnp}^{\text{C}}$ is intractable 
(e.g., due to complex nonlinear dynamics), a control-invariant subset, e.g., 
$\mathcal{Z}^{\text{f,P}}$, is used instead. 
Larger sets yield a larger Parent tube, increasing    
the flexibility of \cmpc\ in exploring its search space.
\end{remark}

As discussed in \cite{LipsitzDownsied}, designing tubes based solely on the Lipschitz constant (which 
in \pmpc\ determines the bound of the uncertainties) 
often leads to very large tubes and excessive conservativeness, 
resulting in overly tightened state constraints.  
While this limits the appeal of conventional TMPC schemes, 
as we demonstrate in Section~\ref{sec:Method}, the Parent-Child MPC architecture 
inverts this shortcoming into an advantage, 
where a deliberately larger tube relaxes the constraints on \cmpc, allowing it to explore 
more dynamic and adaptive control strategies, ultimately resulting in improved system performance.%

\subsubsection{Balancing the computational overhead of \pmpc}
\label{sec:computational_demand_PMPC}

Since the model of \pmpc\ is linear and its nominal optimization problem is convex 
(not necessarily quadratic due to the terminal constraint), \pmpc\ can, 
in view of computational costs, accommodate a larger prediction horizon than \cmpc. 
Larger controllable sets can be reached by sufficiently 
increasing the prediction horizon \cite{Invariant}. To reduce computational burden of \pmpc\ while still enabling  
a large Parent prediction horizon $\Pnp$, we propose two complementary strategies, outlined next.%

\paragraph{Reduced update rate:} 
\pmpc\ does not control the system, so it may be solved less frequently than \cmpc to reduce computational overhead. 
This requires that the updated Parent tube always covers the full Child prediction horizon, i.e., 
the number of time steps between two consecutive control updates by \pmpc\ 
be strictly smaller than $\Pnp - \Cnp$.

\paragraph{Coarser planning sampling:}    
Using larger sampling times, \pmpc\ can extend its horizon 
without adding more optimization variables. 
For stability analysis, the underlying dynamics should use the base sampling time. 
This can be enforced by holding control inputs constant over fixed intervals, 
with the associated variables and constraints removed at implementation.%


\subsection{Child MPC}
\label{sec:C-MPC}

\cmpc\ solves a nonlinear TMPC problem to directly steer the system, 
optimizing the same cost function structure as \pmpc, instead of tracking its state trajectory.   
\pmpc\ shares its nominal control input trajectory and tube with \cmpc, 
which must ensure that its nominal state trajectory remains within the tube 
defined by \pmpc.  
By exploiting its flexibility, \cmpc\ discovers trajectories that enhance optimality 
beyond the conservative solution provided by \pmpc. 
The nominal optimization problem solved by \cmpc\ is formulated by:
\begin{subequations}\label{eq:Child}
\begin{align}
    \begin{split}\label{cost_fun2}
        V\chil&\left(\boldsymbol{x}_k,\tra{\boldsymbol{z}}{k}{k+{\Pnp}}\pare,\tra{\boldsymbol{v}}{k}{k+{\Pnp-1}}\pare\right)= \\
        &\min_{\tra{\boldsymbol{z}}{k}{k+{\Cnp}}\chil,\; 
        \Delta\tra{\boldsymbol{v}}{k}{k+{\Cnp}-1}\chil}
        \Bigg\{\sum^{k + \Cnp - 1}_{i=k}
        \left(
        {\xC_{i|k}}^{\scriptstyle{\top}} Q \; \xC_{i|k}+ {\uC_{i|k}}^{\scriptstyle{\top}} R\; \uC_{i|k}
        \right)+ \\
        &\hspace{25ex} {\xC_{k + \Cnp|k}}^{\scriptstyle{\top}} P\; \xC_{k+\Cnp|k} \Bigg\}
    \end{split}
    \\
       &\text{such that for $i=k,...,k + \Cnp-1$:}
\nonumber
    \\
     \begin{split}\label{initi2}
         \bx_{k}\in  \left\{\xC_k\right\} \oplus \EC 
     \end{split}
    \\
    \begin{split}\label{model2}
        \xC_{i+1|k} = A\xC_{i|k}+g\left(\xC_{i|k}\right)+B\uC_{i|k} 
    \end{split}
    \\
    \begin{split}\label{state2}
        \xC_{j|k}\in \left\{\xP_{j|k} \right\}\oplus \EP \qquad \text{ for }j=k,...,k+\Cnp
    \end{split}
    \\
    \begin{split}\label{input_increment}
        \uC_{i|k}= \uP_{i|k} + \KP \left(\boldsymbol{x}\pare_{i|k}-\xP_{i|k} \right)+\Delta \uC_{i|k} 
    \end{split}
    \\
    \begin{split}
    \label{eq:child_nominal_input}
        \Delta \uC_{i|k}  \in \left(\Vc \ominus \KP \EP \right) \ominus \left\{ \uP_{i|k} \right\}
    \end{split}
\end{align}
\end{subequations}
This optimization is warm-started by the (virtual) solution of \pmpc, i.e., 
by $\tra{\boldsymbol{z}}{k}{k+{\Cnp}}\chil = \tra{\boldsymbol{x}}{k}{k+{\Cnp}}\pare$ and 
$\Delta\tra{\boldsymbol{v}}{k}{k+{\Cnp}-1}\chil = [0, \ldots, 0]^\top$ --- both spanning 
the first $\Cnp$ time steps. 

Constraint \eqref{initi2} ensures that current nominal state $\xC_{k}$ of \cmpc, 
included as a decision variable, is determined  
such that current measured state $\bx_k$ of the system is inside  
the tube $\left\{\xC_k\right\} \oplus \EC$ of \cmpc.
Constraint \eqref{model2} implies that the nominal component of 
\cmpc\ operates based on the original nonlinear dynamic model \eqref{eq:system}, 
omitting the external disturbance $\boldsymbol{w}_k$. Thus, the tube $\EC$ of \cmpc\ 
is designed to account for rejecting these disturbances for the actual system. 
In essence, the Parent-Child MPC architecture allows \cmpc\ 
to avoid designing and deploying explicit terminal sets, a typically challenging task for nonlinear MPC frameworks. 
Instead, a sequence of tightened stage-wise state constraints, 
as given in \eqref{state2}, which ensures that nominal state trajectory of \cmpc\ remains within the 
tube of \pmpc. 
These constraints, when satisfied alongside constraint \eqref{state} within the \pmpc\ nominal 
optimization loop, ensure safe tightening of the \cmpc\ nominal state set, such that 
the actual states remain within the original admissible set $\mathcal{X}$.  
Constraint \eqref{input_increment} defines the nominal \cmpc\ law, 
where the \pmpc\ input \eqref{eq:ancillary} acts as a baseline, and 
\cmpc\ input increment $\Delta\uC_{i|k}$ is treated as a decision variable. 
Constraint \eqref{eq:child_nominal_input} ensures that the resulting 
\cmpc\ nominal inputs remain within the tightened admissible input set $\Vc$.%

In the event that 
$\bx_{k}\notin \left(\mathcal{Z}_{\Cnp}^{\text{C}} \oplus \EC \right)$ holds,  
\cmpc\ invokes \pmpc\ 
to re-solve its optimization problem, initialized at the most recent admissible state,  
to re-define the Parent tube.%

\subsection{Stability and recursive feasibility}

An MPC optimization problem is recursively feasible if, whenever 
a solution exists for an initial state, a feasible solution also exists for all subsequent states 
along the resulting closed-loop trajectory \cite{MPCbook}. 
A given robust positive invariant set $\mathcal{O}$ containing the origin is robustly exponentially stable 
for system \eqref{eq:system}, 
if, for all admissible disturbances $\boldsymbol{w}_k \in \mathcal{W}$, state trajectories starting  
within a certain bounded region converge exponentially towards set $\mathcal{O}$ \cite{LTMPC}. 
The terminal region $\mathcal{Z}^{\text{f,C}}$ defined for \cmpc\ serves as such a 
robust positive invariant set $\mathcal{O}$ within the Parent-Child MPC architecture. 
Thus, once the \cmpc\ nominal trajectory enters set $\mathcal{Z}_{\Cnp}^{\text{C}}$,  
a dual control strategy --- such as the original TMPC formulation in \eqref{eq:TMPC}) --- may be deployed. 
Accordingly, we next establish that the Parent-Child MPC architecture, 
governed by\eqref{eq:Parent}--\eqref{eq:Child}, is guaranteed to exponentially steer the nominal state trajectory 
towards set $\tilde{\mathcal{Z}}_{\Cnp}^{\text{C}}$. 
This, in turn, guarantees that the actual state trajectory reaches the actual terminal set $\mathcal{X}^{\text{f}}$. 
To do so, we should show that the Parent tube, which is by design a robust positive invariant set 
for \pmpc, also remains robustly invariant under the full Parent-Child MPC architecture. 
We accordingly proceed by stating and proving the following two theorems.%

\begin{theorem}
\label{theorem:recursive_feasibility}
    \textbf{Robust invariance of the modeling error set $\EP$ for \cmpc\ nominal state trajectory:} 
    Under the Parent-Child MPC architecture, the modeling error set  $\EP$ --- 
    defining the cross section of the Parent tube --- 
    constitutes a robust positive invariant set for the nominal state 
    trajectory of \cmpc.  
\end{theorem}
\begin{proof}\label{proof1}    
    The stage-wise state constraint \eqref{state2} 
    ensures that the nominal state trajectory of \cmpc\ remains entirely within the Parent tube, 
    making this tube a robust positive invariant set for this trajectory. 
    Designing the modeling error set according to \eqref{eq:EP_stability}, 
    guarantees that the nominal \pmpc\ trajectory is steerable to the origin. 
    Once this occurs, the Parent tube effectively becomes the untranslated modeling error set $\EP$. 
    Consequently, $\EP$ itself serves as a robust positive invariant set for the nominal \cmpc\ trajectory, completing the proof.
    \qed
\end{proof}
For Theorem~\ref{theorem:recursive_feasibility} to be  applicable, 
the nominal \cmpc\ problem should admit at least one feasible trajectory. 
The following theorem establishes the guaranteed existence of such a trajectory.
\begin{theorem}\label{theorem:warm_start}
\textbf{Feasibility of the \pmpc-based warm start for nominal \cmpc:} 
    The solution provided as a warm start for the \cmpc\ nominal problem \eqref{eq:Child} 
    by the \pmpc\ layer, under the ancillary control law \eqref{eq:ancillary}, 
    always yields a feasible solution for nominal \cmpc ---  
    ensuring feasibility for the nominal \cmpc\ problem, whenever the \pmpc\ problem is itself feasible.
\end{theorem}
\begin{proof} 
To prove feasibility, one should show that the warm-started trajectories 
$\tra{\boldsymbol{z}}{k}{k+{\Cnp}}\chil = \tra{\boldsymbol{x}}{k}{k+{\Cnp}}\pare$ and 
$\tra{\boldsymbol{v}}{k}{k+{\Cnp}-1}\chil = \tra{\boldsymbol{u}}{k}{k+{\Cnp}-1}^{\text{P}}$ 
satisfy all constraints of the \cmpc\ problem, i.e., \eqref{initi2}--\eqref{state2} and 
\eqref{eq:child_nominal_input}, with the control law given by \eqref{input_increment}. 
Assuming feasibility for the previous time step $k-1$,  
from \eqref{initi2}, $\bx_{k-1}\in  \left\{\xC_{k-1}\right\} \oplus \EC$ holds. 
The tube of \cmpc\ is robust positive invariant for the actual state trajectory. 
Hence, $ \bx_{k}\in \left\{\xC_{k|k-1}\right\} \oplus \EC$ also holds.
Moreover, based on \eqref{eq:Parent} for the \pmpc\ problem $\xC_{k|k-1} = \boldsymbol{x}_k^{\text{P}}$ holds.  
Therefore, $\bx_{k}\in  \left\{\boldsymbol{x}_k^{\text{P}}\right\} \oplus \EC$ 
holds.
Moreover, with the warm-start input sequence, i.e., 
$\tra{\boldsymbol{v}}{k}{k+{\Cnp}-1}\chil = \tra{\boldsymbol{u}}{k}{k+{\Cnp}-1}^{\text{P}}$,  
constraint \eqref{eq:child_nominal_input} boils down to the \pmpc\ input constraint \eqref{input}, 
which is satisfied by assumption of feasibility of the \pmpc\ problem. 
The dynamics of the warm-started \cmpc\ trajectory aligns with the \pmpc\ dynamics 
\eqref{eq:parent_model}, which is equivalent to 
constraint \eqref{model2} in \cmpc\ when warm-started by $\tra{\boldsymbol{x}}{k}{k+{\Cnp}}\pare$. 
Finally, the tube $\left\{\xP_{i|k} \right\}\oplus \EP$ of \pmpc\ is robust positive invariant, 
i.e., it contains all actual \pmpc\ states --- implying satisfaction of constraint \eqref{state2}.%
    \qed
\end{proof}
Theorems~\ref{theorem:recursive_feasibility} and \ref{theorem:warm_start} jointly 
establish that the Parent-Child MPC architecture guarantees the existence of at 
least one feasible solution that steers system \eqref{eq:system} into a designated terminal set, 
thereby ensuring stabilizability.%
%

\subsection{Generalization to multi-level Parent-Child architecture}

While a single \pmpc\ layer can extend the planning horizon and enhance feasibility, 
it may still fall short in steering the system to the terminal set --- 
particularly in scenarios that demand large horizons and many decision variables due to complex dynamics or tight constraints. 
To address this, a multi-level Parent-Child MPC architecture can be adopted,  
where multiple \pmpc\ systems are stacked in a hierarchical structure. 
In such a recursive structure, each \pmpc\ serves as the Child to a higher-level Parent 
and itself operates as a TMPC augmented with an ancillary control law, as described in Section~\ref{sec:LTMPC}. 
Since each level adheres to the original design assumptions, this architecture remains theoretically valid at all depths.%

As with the bi-level scheme, once the nominal state trajectory of a lower-level 
\pmpc\ enters the untranslated tube of its Parent, 
the top-most Parent layer is discarded, and control is delegated to the next layer. 
For intermediate layers, the required tube is generally larger than that used for the original \cmpc.%

However, adding \pmpc\ layers introduces a critical trade-off: 
While larger Parent tubes facilitate feasibility and robustness, adding more 
layers increases the number of tubes where their Minkowski sum must remain within 
the state admissible set $\mathcal{X}$. 
Therefore, the design must balance warm-start feasibility, computation, and constraint satisfaction across all levels.%

\subsection{Special cases}

Thus far, we have focused on the general implementation of the Parent-Child MPC architecture 
under the assumption that the original system can be stabilized via conventional tube MPC 
using a linear model with bounded modeling error. This section considers simpler, yet common, special cases.%

\subsubsection{Deterministic MPC}

A widely studied variant is the deterministic MPC problem, where the disturbance in \eqref{eq:system} 
is assumed to be zero. In this case, \cmpc\ becomes a standard deterministic MPC and solves the following 
deterministic optimization problem:
\begin{subequations}
\label{deterministic_Child}
\begin{align}
    \begin{split}
        &V\chil\left(\boldsymbol{x}_k\right)= 
        \min_{\tra{\boldsymbol{x}}{k+1}{k+\Cnp}, \Delta\tra{\boldsymbol{u}}{k}{k+\Cnp-1}}
        \Bigg\{\sum^{k + \Cnp - 1}_{i=k}
        \left(
        \boldsymbol{x}^\top_iQ\boldsymbol{x}_i+\boldsymbol{u}^\top_iR\boldsymbol{u}_i
        \right)+ \\
        &\hspace{30ex}
         {\boldsymbol{x}^\top_{k + \Cnp|k}}^{\scriptstyle{\top}} P\; \boldsymbol{x}^\top_{k+\Cnp|k} \Bigg\}
    \end{split}
    \\
       &\text{such that for $i=k,...,k + \Cnp-1$:}
\nonumber
    \\
    \begin{split}
        \boldsymbol{x}_{i+1 | k} = f\left( \boldsymbol{x}_{i|k} , \boldsymbol{u}_{i|k} \right) 
    \end{split}
    \\
    \begin{split}
        \boldsymbol{x}_{i+1|k}\in \left\{\xP_{i + 1|k} \right\}\oplus \EP
    \end{split}
    \\
    \begin{split}
        \boldsymbol{u}_{i|k}=\boldsymbol{u}^{\text{P}}_{i|k} +\Delta \boldsymbol{u}_{i|k} 
    \end{split}
    \\
    \begin{split}
        \boldsymbol{u}_{i|k}  \in \mathcal{U}
    \end{split}
\end{align}
\end{subequations}
The formulation of \pmpc\ remains unchanged, except for constraints \eqref{state} and \eqref{input}, which 
will not involve any disturbance-driven error set (i.e., $\xP_{i+1}\in \mathcal{X} \ominus\EP$ and 
$\uP_i\in \mathcal{U}\ominus K\EP$ for $i=k,...,k + \Pnp-1$). 
Finally, in the second phase, when the actual state enters $\mathcal{Z}_{\Cnp}^{\text{C}}$, 
a deterministic MPC controller replaces the Parent-Child MPC architecture.%

\subsubsection{Linear TMPC}
\label{sec:LTMPC}

Even without disturbances and modeling errors, the MPC problem 
can become computationally challenging when the number of decision variables or 
constraints is large, causing solvers to fails within limited time budgets. 
In such cases, the Parent-Child MPC architecture offers value. 
Since both Parent and Child components share the same system model, 
\pmpc\ can operate with a larger horizon and/or smaller sampling time, 
effectively planning further ahead, especially when the \pmpc\ optimization 
does not need to be solved every iteration. 
Despite the absence of modeling error, a Parent tube with cross section $\EP$ is designed as   
a positive invariant set (rather than a robust one), as well as an ancillary control law  
that establishes the positive invariance of $\EP$. Using solvers for convex or 
quadratic programs ensures computational tractability of \cmpc.%

\section{Case studies}
\label{sec:implemenation}

This section presents two numerical case studies designed to demonstrate the performance 
of the Parent-Child MPC architecture. The first case study highlights its advantage over standard TMPC 
in linear systems, particularly in terms of feasibility and robustness. 
The second case study shows how this architecture can restore feasibility in a nonlinear 
control scenario, where standard TMPC fails due to reachability limitations.
\begin{itemize}
    \item 
    \textbf{Case study 1:} A Parent-Child MPC architecture is implemented for a linear system 
    showcasing its ability to maintain robustness against disturbances while improving feasibility 
    and reducing conservativeness compared to standard TMPC (The code is published in \cite{code1}).
    \item 
    \textbf{Case study 2:} Based on the nonlinear TMPC formulation in \cite{NTMPC}, 
    this case modifies the problem such that the original TMPC becomes infeasible 
    due to a restrictive terminal constraint. The Parent-Child MPC architecture is then 
    introduced to extend the reachability set and to recover feasibility 
    (The code is provided in \cite{code2}).
\end{itemize}
All simulations were conducted on a single computer equipped with an Intel Core i9-13900HX CPU, 16 GB of RAM and an NVIDIA GeForce RTX 4070 GPU.%

\paragraph{\textbf{Case study 1: Linear system}}
We consider a discrete-time double integrator system with additive disturbance, given by:
\begin{equation}
    \boldsymbol{x}_{k+1}=\begin{bmatrix}
        1&1\\ 0&1
    \end{bmatrix}\boldsymbol{x}_{k}+\begin{bmatrix}
        0\\1
    \end{bmatrix}{u}_k+\boldsymbol{w}_k
\end{equation}
This system models a simple kinematic scenario, with the state vector including position and velocity, 
and the control input corresponding to acceleration. 
The control objective is to steer the system from the initial state $\boldsymbol{x}_0 = [2700,0]^{\top}$ 
toward the origin, minimizing cost function:
\begin{equation}\label{eq:cost_linear}
    l(\tra{\boldsymbol{x}}{0}{\Ts},\tra{{u}}{0}{\Ts})=\sum_{i=0}^{\Ts} \left(\boldsymbol{x}^\top_i \boldsymbol{x}_i+u_i^2\right)
\end{equation}
where $\Ts$ is the simulation time. The system is subject to constraints  
$-50\leq \bx[1]\leq10000$, $|\bx[2]| \leq 10000$, $|u|\leq5$, where $[i]$ is used to show element $i$ of a vector. 
The disturbance vector $\boldsymbol{w}_k$ is unknown and bounded $|\boldsymbol{\boldsymbol{w}_k}| \leq 1$ for 
all integer values of $k$.

\begin{figure}
\center
\includegraphics[width=1\columnwidth ]{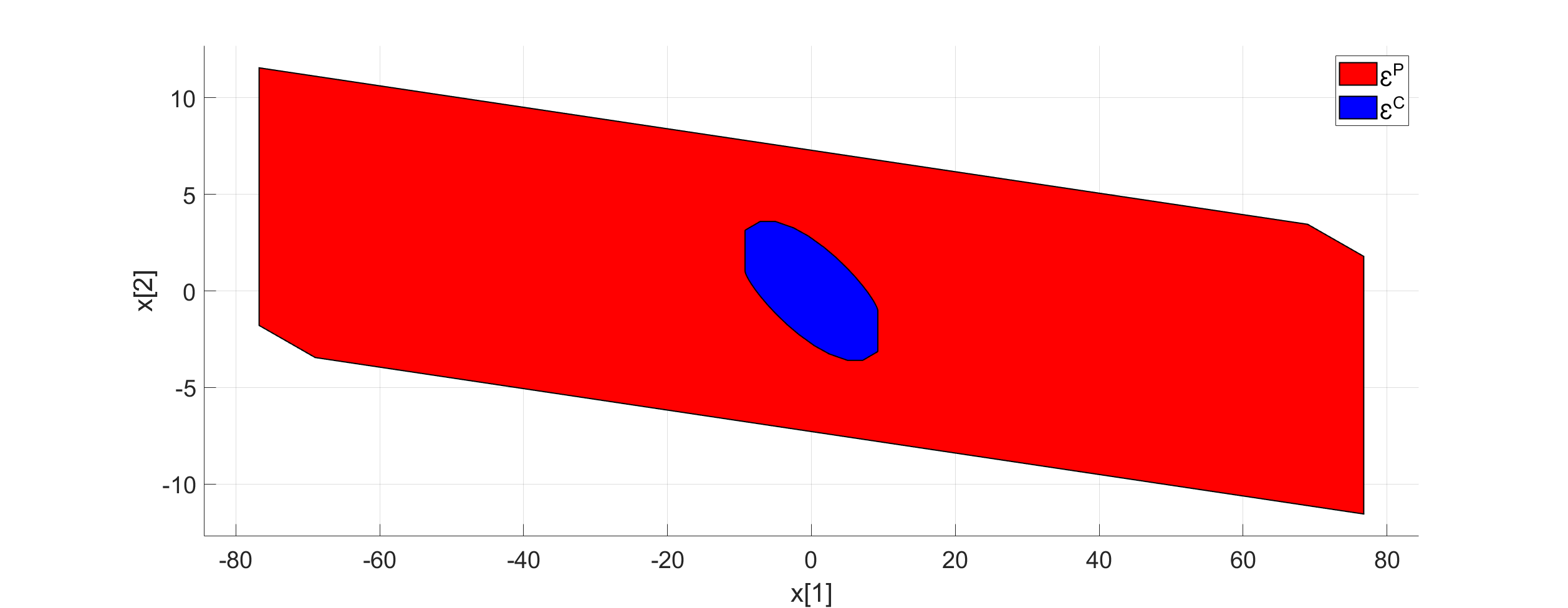}
\caption{Cross section of the disturbance-aware tube used by \cmpc\ (the same as TMPC) and \pmpc.}
\label{fig:tubes}
\end{figure}

The following controllers are implemented for comparison: 

\noindent 
\textbf{TMPC:} 
The simplest controller satisfying the design requirements is a standard linear TMPC, as described in \cite{LTMPC}, 
which outlines conditions for the existence of both the terminal set and the tube. 
We compute the tube via the method in \cite{tube} with the ancillary control law \eqref{eq:ancillary}, 
using the gain $\KP=[0.2054, 0.7835]$, obtained by minimizing the following cost function:
\begin{equation}
    l(\tra{\boldsymbol{z}}{0}{\Ts},\tra{{v}}{0}{\Ts})=\sum_{i=0}^ {\Ts} \left(\boldsymbol{z}^\top_i \boldsymbol{z}_i+10v_i^2\right)
\end{equation}
This cost function penalizes the input more heavily (i.e., by factor $10$) than the original cost in \eqref{eq:cost_linear}. 
This prevents overly aggressive control that will potentially violate input constraints due to excessive 
demand from the ancillary controller.
The corresponding tube cross section is shown in Figure~\ref{fig:tubes}. 
After constraint tightening, the admissible bounds become $-40.77\leq \boldsymbol{z}_k[1] \leq99990.77$,  
$|\boldsymbol{z}_k|\leq9996.41$, and $|v_k|\leq2.86$ for all non-negative integer values $k$.  
We use the MPT3 toolbox \cite{MPT3} to compute the terminal set, which should  
be symmetric about the origin. Since the tightened constraint on $\boldsymbol{z}_k[1]$ 
is asymmetric, it is limited further to its smaller bound, i.e.,  $\left|\boldsymbol{z}_k[1] \right| \leq 40.77$, 
when designing the terminal set. 
The terminal cost matrix in \eqref{cost_fun_TMPC} is set to $P=\begin{bmatrix}
    2.94 & 2.36\\2.36& 4.36 \end{bmatrix}$. 
TMPC uses a horizon of $60$. 

\noindent
\textbf{Bi-level Parent-Child MPC:} 
Prior to reaching $\EP$, nominal \cmpc\ uses stage-wise constraints derived by \pmpc. 
\pmpc\ is designed to guide the system towards the target state $[40,0]^\top$, using the following cost function:
\begin{equation}
    l(\tra{z}{0}{\Ts}\pare,\tra{{v}}{0}{\Ts}\pare)=\sum_{i=0}^{\Ts} \left( ({\xP_i} ^\top-[40,0]) ({\xP_i} -[40,0]^\top)+(v_i^\text{P})^2 \right)
\end{equation}
Targeting $[40,0]^\top$ instead of the origin allows for a larger terminal set $\termPar$ and tube cross section $\EP$, 
delegating the goal of reaching the origin to \cmpc\ (or conventional TMPC in the second control phase). 
Upon entering $\EP$, conventional TMPC is applied with a reduced prediction horizon $10$. 
Since the system is linear and disturbance-free at the Parent level, we can use the MPT3 toolbox to compute both the 
tube cross section and the terminal set. To create the terminal set, we use the Parent 
ancillary control law with gain matrix $\Kf=[0.0057,0.11]$, obtained by solving a feedback design problem with 
a significantly larger input penalty ($\times 10000$) --- allowing  to increase the tube size (see Figure~\ref{fig:tubes}).%

After tightening, the bounds become $36\leq \xP_k[1]\leq99914$, $|\xP_k[2]|\leq9984.86$, and $|v_k\pare|\leq2$ 
for all non-negative integers $k$. 
To reduce the number of decision variables, the same input is repeated every $4$ time steps. 
The Parent prediction horizon $\Pnp$ is set to $120$, which results in $30$ decision variables. 
To further reduce computation, \pmpc\ is re-solved every $6$ time steps.%

\paragraph{\textbf{Results for case study 1}} 
We simulated two $60$-step trajectories, one per controller. 
The Parent-Child MPC architecture required $41$ time steps before \pmpc\ became inactive. 
All optimizations were solved using Matlab's \emph{quadprog} function \cite{quadprog} with the \emph{active set} 
algorithm, which consistently resulted in largest solution speed.%

\begin{figure}
\center
\includegraphics[width=1\columnwidth ]{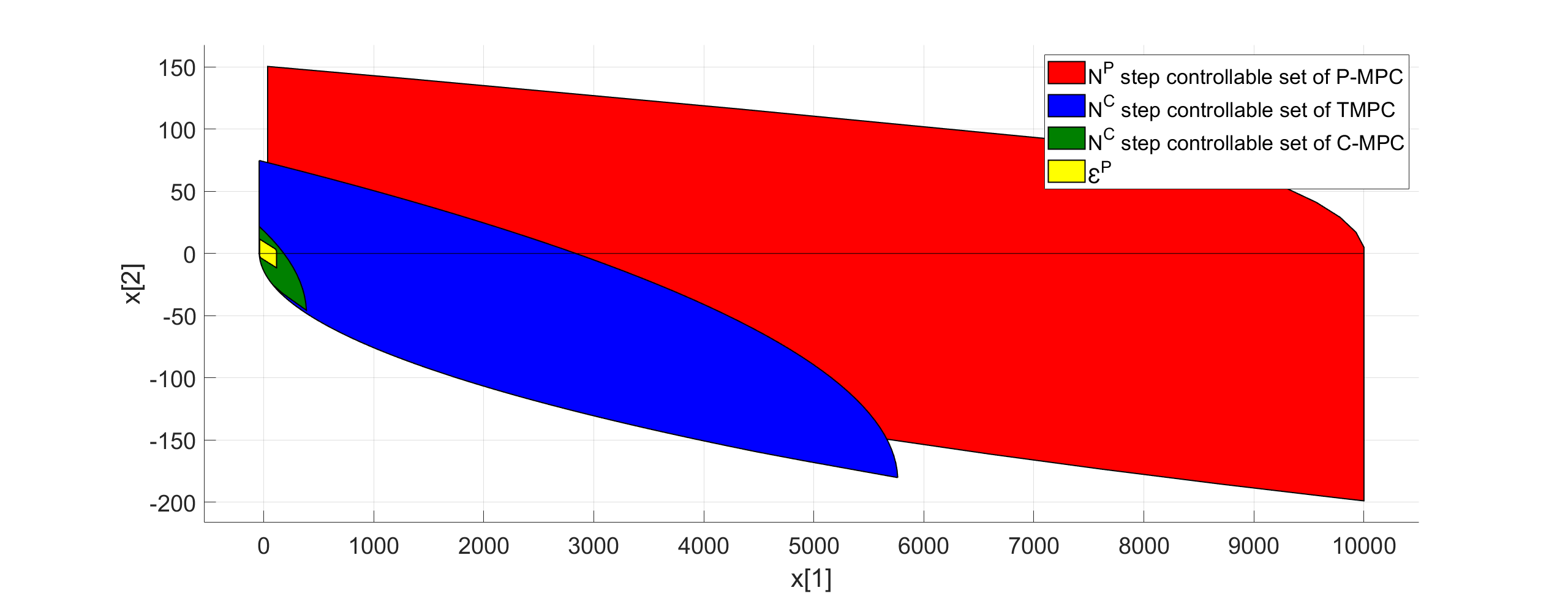}
\caption{Controllable sets and tube cross section for \pmpc.}
\label{fig:reach}
\end{figure}

Figure~\ref{fig:reach} compares three controllable sets, showing that the Parent-Child MPC architecture 
offers the largest set. With zero initial velocity, this architecture successfully solves problems  
considering more than three times the initial distance from the origin that conventional TMPC can handle. 
The figure shows that Parent tube remains within the $\Cnp$-controllable set of \cmpc, 
ensuring stability and recursive feasibility. Conventional TMPC allows enhanced velocity freedom, 
while the tighter acceleration constraints of the Parent-Child MPC architecture results in earlier deceleration.%

\begin{figure*}
\center
\includegraphics[width=2\columnwidth ]{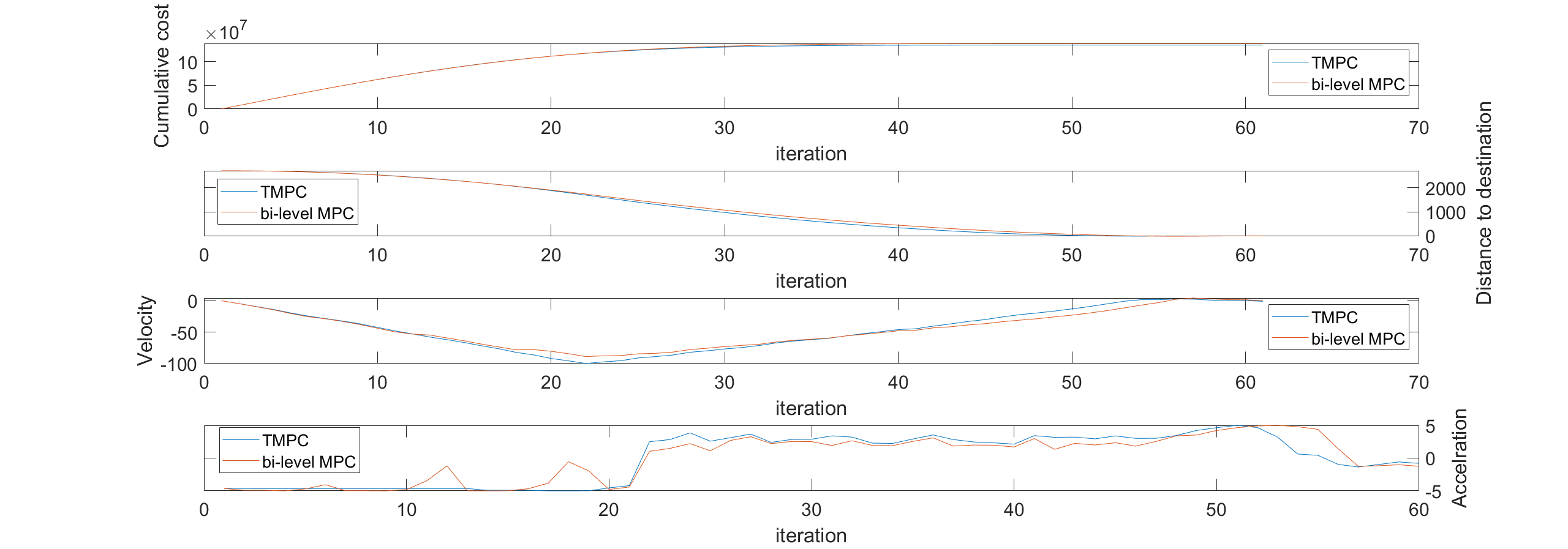}
\caption{Cost, state (including the position and velocity), and control input (i.e., acceleration) over time for TMPC and the Parent-Child MPC architecture.}
\label{fig:tra}
\end{figure*}

Figure~\ref{fig:tra} shows the trajectories of states, control inputs, and cumulative costs:  
TMPC achieves a lower cost, with a maximum advantage of $117.43$ and a final cost $3\%$ lower 
than the Parent-Child MPC architecture, which, however, is significantly faster. 
TMPC takes an average of $8.21$~ms per iteration, compared to $1.10$~ms per iteration for the Parent-Child MPC architecture 
(or $2.73$~ms during \pmpc\ updates).%

\paragraph{\textbf{Case study 2: Nonlinear system}} 
In \cite{NTMPC}, a nonlinear TMPC has been proposed, employing a nonlinear ancillary control law. 
It has been applied to the following system: 
\begin{equation}
    \begin{bmatrix}
        \bx_{k+1}[1]\\ \bx_{k+1}[2]
    \end{bmatrix}=\begin{bmatrix}
        \bx_k[2]\\ \text{sin}\left(\bx_k[1]\right)+u_k
    \end{bmatrix}+\boldsymbol{w}_k
\end{equation}
The state $\boldsymbol{x}_k$ is unconstrained, while the input is limited to $|u_k|\leq0.5$, 
and the disturbance satisfies $|\boldsymbol{w}_k|\leq0.1$. 
The goal was to steer the system from $\boldsymbol{x}_0 = [1, 1]^\top$ towards the origin, 
with a terminal constraint enforcing convergence to zero.%

We consider the same system but tighten the input constraint to $|u_k| \leq 0.3$. 
Under this constraint, the original TMPC (with nominal input limit $|v_k|\leq0.25$) 
cannot reach the terminal set within the specified horizon, resulting in infeasibility.
While extending the horizon potentially resolves this, it significantly increases computational effort 
and, due to non-convexity, can degrade solver performance.%

\begin{figure}
\center
\includegraphics[width=1\columnwidth ]{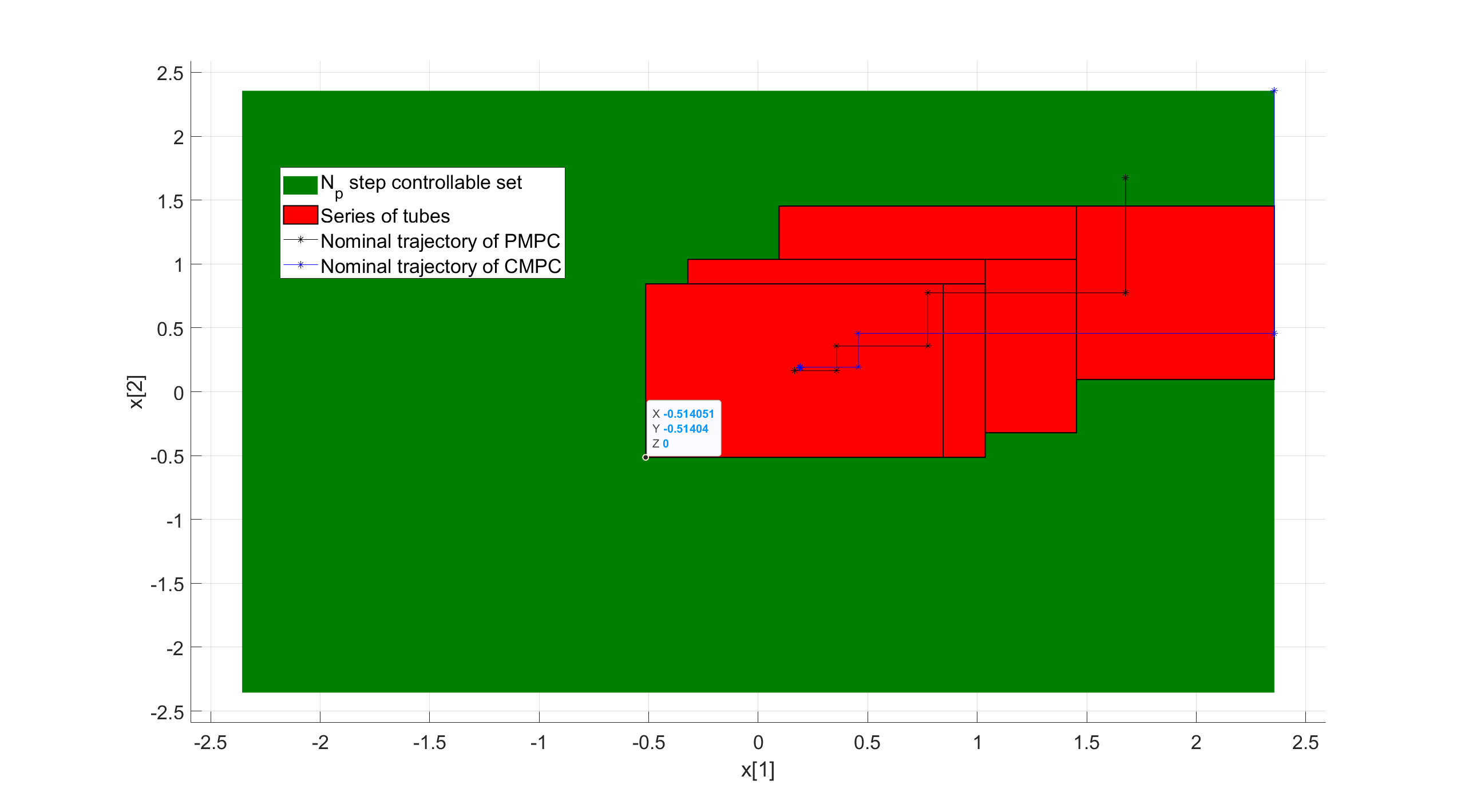}
\caption{The tubes, nominal states of \pmpc, and nominal states of \cmpc.}
\label{fig:init}
\end{figure}

We construct a bi-level Parent-Child MPC architecture using the original TMPC as \cmpc\  
and the following linearized model for \pmpc:
\begin{equation}
    \xC_{k+1}=\begin{bmatrix}
        0&1\\ 0.46 &0
    \end{bmatrix}\xC_k+\begin{bmatrix}
        0\\1
    \end{bmatrix}v_k\chil+\begin{bmatrix}
        0\\g\left(\xC_k\right)
    \end{bmatrix}
\end{equation}
where $g\left(\boldsymbol{x}_k\right) = \sin \left(\bx_k[1]\right) - 0.46 \bx_k[1]$. 
To ensure validity, we restrict the state to $\left|\xC_k\right| \leq 4\pi/3$. 
These constraints do not reduce the feasibility compared to the original TMPC, 
which was already infeasible outside this domain. 
The function $g(\cdot)$ is Lipschitz and bounded, allowing us to treat it as 
a bounded modeling error for \pmpc, i.e., $\wP_k[1] = 0$ and $\left|\wP_k[2]\right| \leq 0.4$.%

A linear ancillary control law is designed using \eqref{eq:ancillary} with $\KP = [0.051, 0]$, 
optimized for $Q = I_{2\times 2}$ and $R = 20$. 
This yields a tube $\mathcal{E}$, such that $|\eP_k| \leq 0.68$. 
Consequently, the nominal input must satisfy $\left|v^\text{P}_k\right| \leq 0.22$. 
The terminal set is defined as $\left|\xP_{k+ \Pnp}\right| \leq 0.68$, with a horizon of $\Pnp = 10$, 
and a terminal cost given by ${\xP_{k+N^\text{P}}}^\top\text{diag}(1.15,2.15)\xP_{k+ \Pnp}$.   
For \cmpc, we reuse the controller from \cite{NTMPC}, replacing the terminal constraint with time-varying constraints from \pmpc, 
as described in \eqref{eq:Child}. 
Optimization problems are solved using \emph{quadprog} or \emph{fmincon} \cite{fmincon} 
with the SQP algorithm for nonlinear problems.%


\begin{figure}
\center
\includegraphics[width=1\columnwidth ]{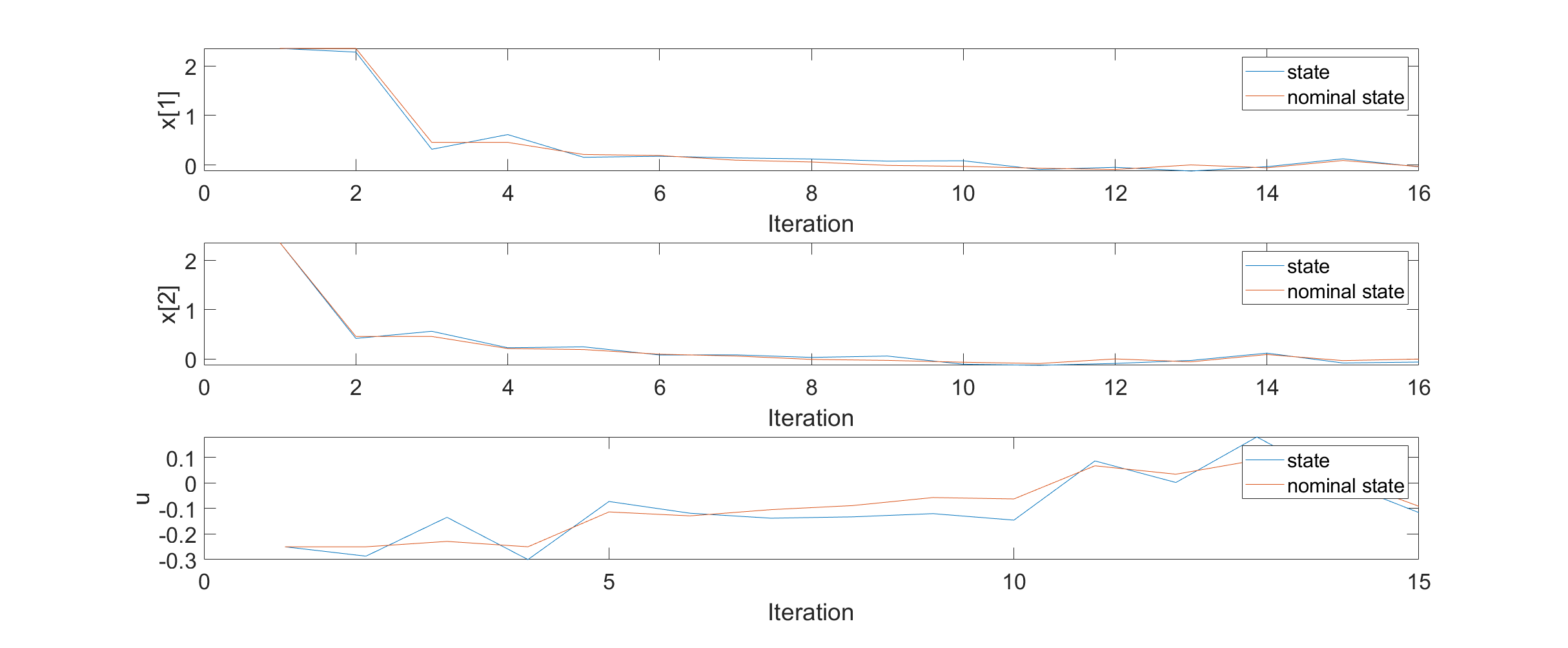}
\caption{Trajectory of nominal and actual states and inputs of \cmpc.}
\label{fig:sim}
\end{figure}

\begin{figure}
\center
\includegraphics[width=1\columnwidth ]{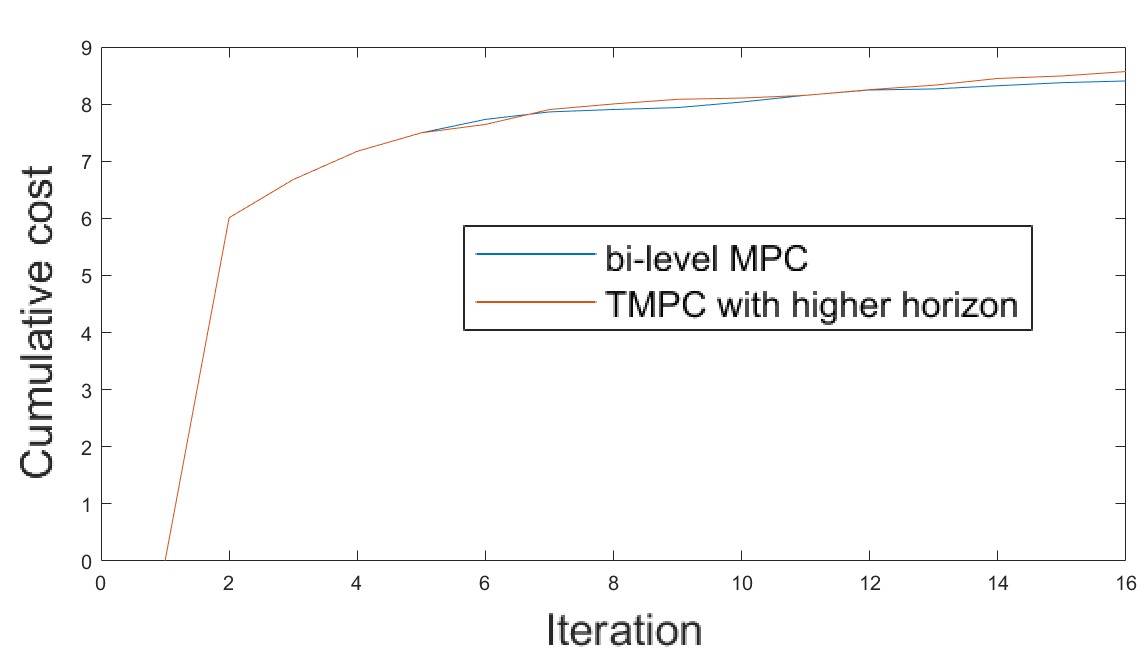}
\caption{Cumulative cost of the Parent-child MPC architecture and original TMPC with a larger horizon.}\label{fig:cost}
\end{figure}

\paragraph{\textbf{Results for case study 2}} 
Both the Parent-Child MPC architecture and conventional TMPC with an extended horizon 
were implemented to steer the system from $\bx_0=[\pi4/3,\pi4/3]^\top$ towards the origin.  
Figure~\ref{fig:init} shows the tubes generated by \pmpc\ and the 
nominal trajectory determined by \cmpc, demonstrated that \cmpc\ 
successfully finds a feasible trajectory. 
Simulation results for the Parent-Child MPC architecture are presented in Figure~\ref{fig:sim}. 
Most notably, Figure~\ref{fig:cost} compares the cumulative costs of both controllers, 
showing that they achieve comparable performance under the same cost function.
Although both controllers yield similar optimal performance, their computational costs differ significantly. 
TMPC with an extended horizon requires an average of $7.2$~ms per solve, 
while \cmpc\ takes $3.~$~8ms and \pmpc\ just $0.5$~ms on average.%

\section{Conclusions and future work}\label{sec:conclusions}

This paper introduced a novel architecture, called the Parent-Child Model Predictive Control (MPC). 
The Parent MPC (\pmpc) layer generates long-term plans using a simplified model and/or smaller sampling 
times for computational efficiency and to ensure stability, while the Child MPC (\cmpc) layer 
refines the plan in real time to enhance performance. 
We provided formal guarantees for stability and recursive feasibility of the Parent-Child MPC architecture.%

The primary advantage of this architecture lies in its ability to reduce computational cost, 
while planning further ahead. Although \cmpc\ solves nonlinear, non-convex problems, 
the warm start from \pmpc\ improves the solver efficiency.%

We demonstrated the benefits of the Parent-Child MPC architecture through two case studies. 
The first showed that this architecture reduces the computation time and expands the controllable set in linear systems. 
The second illustrated that --- thanks to short-term corrections of \cmpc\ --- 
the architecture remains effective for systems with 
Lipschitz nonlinearity, even though \pmpc\ uses an inaccurate model. 

Future research should focus on applying this architecture to more complex controllers, 
such as state-dependent dynamic tube-based MPC \cite{SDDTMPC}, and to real-world systems, 
e.g., quadrotors \cite{quadrotor}. Evaluating and comparing tube-generation algorithms 
to enhance feasibility and performance in more demanding applications are further topics for future research.%





\begin{thebibliography}{10}

\bibitem{MPCbook}
J.~Rawlings, D.~Mayne, and M.~Diehl, {\em Model Predictive Control: Theory, Computation, and Design}.
\newblock Nob Hill Publishing, 2017.

\bibitem{passive}
P.~Falugi, ``Model predictive control: a passive scheme,'' {\em IFAC Proceedings Volumes}, vol.~47, no.~3, pp.~1017--1022, 2014.
\newblock 19th IFAC World Congress.

\bibitem{stateCostModifided}
W.-H. Chen and Y.~Yan, ``New stability theory of model predictive control: modified stage cost approach,'' {\em International Journal of Systems Science}, vol.~56, no.~4, pp.~808--826, 2025.

\bibitem{Invariant}
E.~Kerrigan and J.~Maciejowski, ``Invariant sets for constrained nonlinear discrete-time systems with application to feasibility in model predictive control,'' in {\em Proceedings of the 39th IEEE Conference on Decision and Control (Cat. No.00CH37187)}, vol.~5, pp.~4951--4956 vol.5, 2000.

\bibitem{extendedHorizon}
S.~Liu and J.~Liu, ``Economic model predictive control with extended horizon,'' {\em Automatica}, vol.~73, pp.~180--192, 2016.

\bibitem{quadrotor}
M.~Chipofya, D.~Lee, and K.~Chong, ``Trajectory tracking and stabilization of a quadrotor using model predictive control of laguerre functions,'' {\em Abstract and Applied Analysis}, vol.~2015, 02 2015.

\bibitem{mobileRobot}
Z.~Sun, L.~Dai, K.~Liu, Y.~Xia, and K.~H. Johansson, ``Robust {MPC} for tracking constrained unicycle robots with additive disturbances,'' {\em Automatica}, vol.~90, pp.~172--184, 2018.

\bibitem{sampling}
D.~Zhou and K.~V. Ling, ``The effect of sample/hold time on initial feasible set in model predictive control design,'' in {\em 2013 Australian Control Conference}, pp.~212--217, 2013.

\bibitem{overviewMPC}
K.~Holkar, K.~Wagh, and L.~Waghmare, ``An overview of model predictive control,'' {\em International Journal of Control and Automation}, vol.~3, 01 2011.

\bibitem{DualTerminal}
D.~He, ``Dual-mode nonlinear {MPC} via terminal control laws with free-parameters,'' {\em IEEE/CAA Journal of Automatica Sinica}, vol.~4, no.~3, pp.~526--533, 2017.

\bibitem{Hier}
M.~Brdys, M.~Grochowski, T.~Gminski, K.~Konarczak, and M.~Drewa, ``Hierarchical predictive control of integrated wastewater treatment systems,'' {\em Control Engineering Practice}, vol.~16, no.~6, pp.~751--767, 2008.
\newblock Special Section on Large Scale Systems.

\bibitem{Hier2}
G.~S. van~de Weg, H.~L. Vu, A.~Hegyi, and S.~P. Hoogendoorn, ``A hierarchical control framework for coordination of intersection signal timings in all traffic regimes,'' {\em IEEE Transactions on Intelligent Transportation Systems}, vol.~20, no.~5, pp.~1815--1827, 2019.

\bibitem{Anahita}
A.~Jamshidnejad, D.~Sun, A.~Ferrara, and B.~{De Schutter}, ``A novel bi-level temporally-distributed {MPC} approach: An application to green urban mobility,'' {\em Transportation Research Part C: Emerging Technologies}, vol.~156, p.~104334, 2023.

\bibitem{FLMPC}
F.~Surma and A.~Jamshidnejad, ``Fuzzy-logic-based model predictive control: A paradigm integrating optimal and common-sense decision making,'' {\em ArXiv}, 2025.

\bibitem{setComputation}
F.~Tahir, ``Efficient computation of robust positively invariant sets with linear state-feedback gain as a variable of optimization,'' in {\em 2010 7th International Conference on Electrical Engineering Computing Science and Automatic Control}, pp.~199--204, 2010.

\bibitem{LTMPC}
S.~Raković and D.~Mayne, ``A simple tube controller for efficient robust model predictive control of constrained linear discrete time systems subject to bounded disturbances,'' {\em IFAC Proceedings Volumes}, vol.~38, no.~1, pp.~241--246, 2005.
\newblock 16th IFAC World Congress.

\bibitem{epuck}
P.~Gonçalves, P.~Torres, C.~Alves, F.~Mondada, M.~Bonani, X.~Raemy, J.~Pugh, C.~Cianci, A.~Klaptocz, S.~Magnenat, J.-C. Zufferey, D.~Floreano, and A.~Martinoli, ``The e-puck, a robot designed for education in engineering,'' {\em Proceedings of the 9th Conference on Autonomous Robot Systems and Competitions}, vol.~1, 01 2009.

\bibitem{SDDTMPC}
F.~Surma and A.~Jamshidnejad, ``State-dependent dynamic tube {MPC}: A novel tube {MPC} method with a fuzzy model of disturbances,'' {\em International Journal of Robust and Nonlinear Control}, vol.~35, no.~4, pp.~1319--1354, 2025.

\bibitem{NTMPC}
D.~Q. Mayne, E.~C. Kerrigan, E.~J. van Wyk, and P.~Falugi, ``Tube-based robust nonlinear model predictive control,'' {\em International Journal of Robust and Nonlinear Control}, vol.~21, no.~11, pp.~1341--1353, 2011.

\bibitem{LipsitzImplmented}
H.~E.~T. Yiqi~Gao, Andrew~Gray and F.~Borrelli, ``A tube-based robust nonlinear predictive control approach to semiautonomous ground vehicles,'' {\em Vehicle System Dynamics}, vol.~52, no.~6, pp.~802--823, 2014.

\bibitem{Lipschitz}
S.~Yu, H.~Chen, and F.~Allgöwer, ``Tube {MPC} scheme based on robust control invariant set with application to lipschitz nonlinear systems,'' in {\em 2011 50th IEEE Conference on Decision and Control and European Control Conference}, pp.~2650--2655, 2011.

\bibitem{tube}
S.~Rakovic, E.~Kerrigan, K.~Kouramas, and D.~Mayne, ``Invariant approximations of the minimal robust positively invariant set,'' {\em IEEE Transactions on Automatic Control}, vol.~50, no.~3, pp.~406--410, 2005.

\bibitem{LipsitzDownsied}
S.~Singh, A.~Majumdar, J.-J. Slotine, and M.~Pavone, ``Robust online motion planning via contraction theory and convex optimization,'' in {\em 2017 IEEE International Conference on Robotics and Automation (ICRA)}, pp.~5883--5890, 2017.

\bibitem{code1}
F.~Surma, ``Code to test the implementation of linear bi-level parent-child mpc,'' 2025.

\bibitem{code2}
F.~Surma, ``Code to test the implementation of nonlinear bi-level parent-child mpc,'' 2025.

\bibitem{MPT3}
M.~Herceg, M.~Kvasnica, C.~Jones, and M.~Morari, ``{Multi-Parametric Toolbox 3.0},'' in {\em Proc.~of the European Control Conference}, (Z\"urich, Switzerland), pp.~502--510, July 17--19 2013.
\newblock \url{http://control.ee.ethz.ch/~mpt}.

\bibitem{quadprog}
Mathworks, ``quadprog.'' \url{https://www.mathworks.com/help/optim/ug/quadprog.html}.

\bibitem{fmincon}
Mathworks, ``fmincon.'' \url{https://www.mathworks.com/help/optim/ug/fmincon.html}.

\end{thebibliography}

\end{document}